\documentclass[11pt,leqno]{article}
\usepackage{amsmath, amsfonts, amssymb, amsthm}

\theoremstyle{definition}
\newtheorem{theorem}{Theorem}[section]

\newtheorem{proposition}[theorem]{Proposition}

\newtheorem{remark}[theorem]{Remark}

\newcommand{\be}{\begin{equation}}
\newcommand{\ee}{\end{equation}}
\newcommand{\beq}{\begin{equation*}}
\newcommand{\eeq}{\end{equation*}}
\newcommand{\bq}{\begin{eqnarray}}
\newcommand{\eq}{\end{eqnarray}}
\newcommand{\bqn}{\begin{eqnarray*}}
\newcommand{\eqn}{\end{eqnarray*}}

\newcommand\ceil[1]{\lceil#1\rceil}
\newcommand\floor[1]{\lfloor#1\rfloor}

\newcommand{\twopiece}[4]{\left\{\begin{array}{ll}{#1}&{#2}\\{#3}&{#4}\end{array}\right.}
\newcommand{\threepiece}[6]{\left\{\begin{array}{ll}{#1}&{#2}\\{#3}&{#4}\\{#5}&{#6}\end{array}\right.}

\newcommand{\eps}{\varepsilon}

\newcommand{\GG}{\mathbb G}
\newcommand{\RR}{\mathbb R}
\newcommand{\NN}{\mathbb N}

\newcommand{\cD}{{\mathcal D}}

\newcommand{\cF}{{\mathcal F}}
\newcommand{\cG}{{\mathcal G}}

\newcommand{\cS}{{\mathcal S}}

\begin{document}

\title{Entropy Estimate For High Dimensional Monotonic Functions}
\author{{\sc Fuchang Gao}\thanks{Corresponding author. Department of Mathematics, P.O.
Box 441103, University of Idaho, Moscow, ID 83844-1103. Email:
fuchang@uidaho.edu. Phone: 1-208-885-5274. Fax: 1-208-885-5843.}\,
\thanks{Supported in part by NSF Grant DMS-0405855.}
\\ {\it Department of Mathematics}\\{\it University of Idaho} \and
 {\sc Jon A. Wellner}\thanks{Supported in part by NSF
Grant DMS-0503822.}
\\{\it Department of Statistics} \\{\it University of Washington}}
\date{\today}
\maketitle
\begin{abstract}
We establish upper and lower bounds for the metric entropy and bracketing entropy of the class
of $d$-dimensional bounded monotonic functions under $L^p$ norms. It
is interesting to see that both the metric entropy and bracketing
entropy have different behaviors for $p<d/(d-1)$ and $p>d/(d-1)$. We
apply the new bounds for bracketing entropy to establish a global
rate of convergence of the MLE of a $d$-dimensional monotone
density.
\end{abstract}

\bigskip


\newpage

\section{Introduction}
Shape constrained functions appear very commonly in
nonparametric estimation in statistics via renewal theory and
mixing of uniform distributions.
 A class of multivariate functions of interests in
applications is the class of ``block-decreasing'' densities; see
e.g. Polonik \cite{polonik:95}, \cite{polonik:98}, and Biau and
Devroye \cite{biau-devroye:03}. It consists of bounded densities on
$\RR^d$ that are decreasing in each variable. We denote by $\cF_d$
the collection of non-negative functions on $[0,1]^d$ which are
bounded by 1, and monotonic in each variable, that is, monotonic
along any line that is parallel to an axis. As is well known, the
rate of convergence of nonparametric estimators such as the Maximum
Likelihood Estimator (MLE) is determined by the metric entropy and
bracketing entropy bounds for an appropriate related class of
 functions; see the definitions below.

In this paper, we provide upper and lower bounds for the entropy
$\log N(\eps,\cF_d,\|\cdot\|_p)$ and the bracketing entropy $\log
N_{[\,]}(\eps,\cF_d,\|\cdot\|_p)$, where $N(\eps,\cF_d,\|\cdot\|_p)$
and $N_{[\,]}(\eps,\cF_d,\|\cdot\|_p)$ are defined as follows:
$$
 N(\eps, \cF_d,\|\cdot\|_p)
 :=\min \left\{m: \exists f_1,f_2, \dots, f_m \mbox { s.t. } \cF_d\subset \bigcup_{k=1}^m B_p(f_k,\eps )\right\}
$$
where $B_p(f_k,\eps)=\{f\in \cF_d: \|f-f_k\|_p\le \eps\}$, and
$$
N_{[\,]}(\eps, \cF_d,\|\cdot\|_p):=\min \left\{m: \exists
\underline{f}_1,\overline{f}_1, \dots,
\underline{f}_m,\overline{f}_m  \mbox { s.t. }
\|\overline{f}_k-\underline{f}_k\|_p\le \eps, \cF_d\subset
\bigcup_{k=1}^m[\underline{f}_k,\overline{f}_k]\right\},
$$
where
$$
[\underline{f}_k,\overline{f}_k]=\left\{g\in \cF_d:
\underline{f}_k\le g\le \overline{f}_k\right\}.
$$
The new bracketing entropy bounds have implications for the rate of
convergence of the Maximum Likelihood Estimator of a ``block
decreasing'' density as will be shown in section 5.

Our main result is the following
\begin{theorem}
\label{BracketingEntropyBoundTheorem}
For $p\ge 1$, there exist constants $c_1$ and $c_2$ depending only
on $p$, such that if $(d-1)p\ne d$, then
$$
{c_1\eps^{-\alpha}}\le \log N(\eps, \cF_d,\|\cdot\|_p )\le \log
N_{[\,]}(\eps, \cF_d,\|\cdot\|_p )\le {c_2\eps^{-\alpha}},$$ where
$\alpha=\max\{d,(d-1)p\}$. If $(d-1)p=d$, then
\begin{eqnarray}
{c_1\eps^{-d}}\le \log N(\eps, \cF_d,\|\cdot\|_p )\le \log
N_{[\,]}(\eps, \cF_d,\|\cdot\|_p )\le {c_2\eps^{-d}(\log
1/\eps)^{1+d/p}}.
\label{BoundaryCaseBoundsMainTheorem}
\end{eqnarray}
\end{theorem}

\begin{remark}We believe that in the critical case $(d-1)p=d$, the
logarithmic factor in the upper bound in (\ref{BoundaryCaseBoundsMainTheorem})
 is not needed, and prove in
Theorem 4.1 that this is indeed so for regular entropy under the $L^p$
norm, provided $(d,p)\ne (2,2)$.
\end{remark}

It should be pointed out that when $d=1$, $\cF_d$ is just the class
of probability distribution functions, and the entropies are known
to be of the order $\eps^{-1}$; see e.g. \cite{vdv-w:96}, Theorem
2.75, page 159. So, in some sense, the results in this paper
generalize the known results for $d=1$. It should also be noted that
when $d>1$, $\cF_d$ is a much larger class than that of
$d$-dimensional probability distributions. Indeed, Blei, Gao and Li
\cite{BGL:05} recently proved that under the $L^2$ norm, the metric
entropy of the class $\cD_d$ of $d$-dimensional probability
distributions satisfies
$$
c_1\varepsilon^{-1}[\log(1/\varepsilon)]^{d-1/2}(\log\log(1/\eps))^{-1/2}\le\log
N(\varepsilon, \cD_d,\|\cdot\|_2)\le
c_2\varepsilon^{-1}[\log(1/\varepsilon)]^{d-1/2}.
$$
for $d>2$, and
$$
c_1\varepsilon^{-1}[\log(1/\varepsilon)]^{3/2}\le\log N(\varepsilon,
\cD_d,\|\cdot\|_2)\le
c_2\varepsilon^{-1}[\log(1/\varepsilon)]^{3/2}.
$$
for $d=2$.

The paper is organized as follows. First, we prove the lower bound
for regular entropy by constructing a well-separated set using a
combinatorial argument. Next, we obtain the upper bound for
bracketing entropy using a constructive proof, revealing the
difference of entropy growth between the cases $p<d/(d-1)$ and
$p>d/(d-1)$. Then we turn to the critical case $p=d/(d-1)$, and use
the result for the case $p=1$ and the metric entropy estimate of
convex hulls to remove the extra logarithmic factor in the upper
bound for the regular entropy. Finally, we apply the bracketing
entropy estimate to establish a global rate of convergence of the
MLE of a $d$-dimensional ``block-decreasing'' density.

\section{Lower bound}
In this section, we obtain the lower bound estimate, namely
\begin{proposition} For $p\ge 1$, there exists a constant $c_1>0$ such
that
$$
\log N(\eps, \cF_d,\|\cdot\|_p)\ge {c_1\eps^{-\alpha}},
$$
where $\alpha=\max\{d,(d-1)p\}$.
\end{proposition}
\begin{proof}
For convenience, we assume $\eps=2^{-n}$ for some positive integer
$n$. We divide $[0,1]^d$ into $\eps^{-d}$ small cubes of side-length
$\eps$. Define $g$ on $[0,1]^d$, such that on each open cube
$\prod_{i=1}^d(k_i\eps,k_i\eps+\eps)$, $0\le k_i<2^n$, $1\le i\le
d$,
$$
g(x)=\frac{(k_1+k_2+\cdots+k_d+1)\eps}{3d}\pm \frac{\eps}{6d}.
$$
Clearly, there are $2^{\eps^{-d}}$ different ways to define $g$,
and each can be extended to a function in $\cF_d$. Let $\cG_d$ be
the collection of these extended functions.

For each $g\in \cG_d$ define
$$
 B(g)=\{h\in \cG_d: \mbox{ there are
at most $2^{-4}{\eps^{-d}}$ open cubes on which $g\ne h$}\}.
$$
Since ${m \choose l} \le (me/l)^l$ and $(16e)^{1/16} \le 2^{1/2}$,
it is easy to check that $B(g)$ contains no more than
${{\eps^{-d}}\choose{2^{-4}{\eps^{-d}}}}\le 2^{\eps^{-d}/2}$
elements. Thus, we can find $N=2^{\eps^{-d}/2}$ functions $g_1$,
$g_2$, ..., $g_N$, such that if $i\ne j$, then $B(g_i)$ and $B(g_j)$
are disjoint. Clearly
$$
\|g_i-g_j\|_1\ge \frac{\eps}{3d}\cdot
\frac{1}{2^4}=\frac{\eps}{48d}.
$$
Hence, $N((48d)^{-1}\eps, \cF_d,\|\cdot\|_1 )\ge 2^{\eps^{-d}/2}$,
which implies
$$
N(\eps, \cF_d,\|\cdot\|_p,)\ge N(\eps, \cF_d,\|\cdot\|_1,)\ge
e^{c_1\eps^{-d}}
$$
for some constant $c_1>0$ and all $p\ge 1$.

When $p>d/(d-1)$, this lower bound is not sharp. In order to improve
it, we will construct a different well-separated subset. We define
$q(x)$ on $[0,1]^d$ as follows: on each open cube
$\prod_{i=1}^d(k_i\eps,k_i\eps+\eps)^d$ that satisfies
$k_1+k_2+\cdots+k_d=\eps^{-1}$, $k_1, k_2,..., k_d\ge 0$, we define
$q(x)=\frac12\pm \frac12$. Clearly, $q(x)$ can be extended to a
function in $\cF_d$. Now, because there are $c\eps^{1-d}$ qualified
cubes, where $c$ is a constant depending only on $d$, there are
$2^{c\eps^{1-d}}$ different functions $q(x)$. The same combinatorial
argument as the one given above shows that there are at least
$m=2^{c\eps^{1-d}/2}$ functions $q_1$, $q_2$, ..., $q_m$, such that
$|q_i-q_j|=1$ on at least $c\eps^{1-d}/2^4$ cubes, $i\ne j$. Thus,
$$
\|q_i-q_j\|_p\ge \left(\frac{c\eps}{2^4}\right)^{1/p}.
$$
This implies that
$$
N(\left(c2^{-4}\eps \right)^{1/p}, \cF_d,\|\cdot\|_p)\ge
2^{c\eps^{1-d}/2},$$ which further implies
$$
N(\eps, \cF_d,\|\cdot\|_p)\ge e^{c_1\eps^{-(d-1)p}},
$$
for some constant $c_1>0$ when $p>d/(d-1)$.
\end{proof}
\section{Upper bound}
In this section, we obtain an upper bound through a constructive
proof. We will prove
\begin{proposition} For $p\ge 1$, $p\ne d/(d-1)$, there exists a constant $c_2>0$ such
that
$$
\log N_{[\,]}(\eps, \cF_d,\|\cdot\|_p)\le {c_2\eps^{-\alpha}},
$$
where $\alpha=\max\{d,(d-1)p\}$. For $p=d/(d-1)$, there exists a
constant $c_2>0$ such that
$$
\log N_{[\,]}(\eps, \cF_d,\|\cdot\|_p)\le {c_2\eps^{-d}(\log
1/\eps)^{1+d/p}}.
$$
\end{proposition}
\subsection{Construction}
For convenience, we introduce the notion
$$
\omega(f,I)=\sup\{f(t): t\in I\}-\inf\{f(t):t\in I\},$$ where $I$ is
any subset of $[0,1]^d$.

If $p=1$, we choose $K=2^d$; otherwise, we choose $K=2^\beta$ where
$\beta=\frac{1}{2}[d-1+1/(p-1)]$. For any given $\eps=2^{-n}$, $n\in
\NN$, let $l$ be the integer satisfying $K^{-l}\le \eps<K^{-l+1}$.

For each $f\in \cF_d$, we construct $\underline{f}$ and
$\overline{f}$ as follows. First, we partition $[0,1)^d$ into
$\eps^{-d}$ cubes of side-length $\eps$. (All the cubes are of the
form $\prod_{i=1}^d[a_i,b_i)$.) A cube $I_0$ of side-length $\eps$
is selected if $\omega(f,I_0)\le K\eps$. For each cube that is not
selected, we partition it into $2^d$ cubes of equal size. In
general, suppose we have a cube $I_i$ of side-length $2^{-i}\eps$.
If $\omega(f,I_i)\le K^{i+1}\eps$, we select the cube; otherwise we
partition the cube into $2^d$ smaller cubes. This process continues
until $i=l$. In this case, we always select the cube. Clearly, each
point in $[0,1)^d$ uniquely belongs to one of the selected cubes.

On each selected cube $I$ of side-length $2^{-i}\eps$, $0\le i<l$,
we define
$$
\underline{f} =K^{i+1}\eps\left\lfloor\frac{\inf_{x\in
I}f(x)}{K^{i+1}\eps}\right\rfloor,\,\,\,\,\,\,\,
\overline{f}=K^{i+1}\eps\left\lceil{\frac{\sup_{x\in
I}f(x)}{K^{i+1}\eps}}\right\rceil.
$$
On each selected cube of side-length $2^{-l}\eps$ and on
$[0,1]^d\setminus [0,1)^d$, we define $\overline{f}=1$ and
$\underline{f}=0$. Clearly, $\underline{f}\le f\le \overline{f}$.

Let $\overline{\cS}=\{\overline{f}: f\in \cF_d\}$, and
$\underline{\cS}=\{\underline{f}: f\in \cF_d\}$. We will estimate
$\|\overline{f}-\underline{f}\|_p$, and the cardinalities
$|\underline{\cS}|$ and $|\overline{\cS}|$ of $\underline{\cS}$ and
$\overline{\cS}$ respectively.

\subsection{Bound for $\|\overline{f}-\underline{f}\|_p$}
For each $i\in \NN$, let $U_i$ be the union of the selected cubes of
side-length $2^{-i}\eps$. We first bound the measure of $U_i$.

Let $s_i$ be the number of cubes of side-length $2^{-i}\eps$ that
have been selected, and $n_i$ be the number of cubes of side-length
$2^{-i}\eps$ that have not been selected. Clearly, by the
construction of $\underline{f}$ and $\overline{f}$, we have
$s_i+n_i=2^dn_{i-1}$. In particular, $s_i\le 2^dn_{i-1}$.

Now we try to estimate $n_{i-1}$ for $i\ge 1$. If a cube
$I=\prod_{j=1}^d[a_j,b_j)$ of side-length $2^{-i+1}\eps$ is not
selected, then $\omega(f,I)>K^{i}\eps$. By the monotonicity of $f$
along each variable, there exists $1\le j\le d$, such that on the
edge $\overline{A_{j-1}A_{j}}$, we have
$\omega(f,\overline{A_{j-1}A_{j}})>K^{i}\eps/d$, where
$$
A_{j}=(b_1,...,b_{j},a_{j+1},...,a_d) .
$$
 Thus for $n_{i-1}$
cubes of side-length $2^{-i+1}\eps$, there are $n_{i-1}$ disjoint
edges on which $\omega(f,\cdot)>K^{i}\eps/d$. From these edges,
there are at least $\ceil{n_{i-1}/d}$ edges that are parallel.
Furthermore from these parallel edges, there are at least
$\ceil{n_{i-1}(2^{-i+1}\eps)^{d-1}/d}$ disjoint edges that lie on
the same line segment $[0,1]$ that is parallel to one of the axes.
Because $f$ is monotonic along this line segment, and the value
change is at most $1$, we have
$$
\ceil{n_{i-1}(2^{-i+1}\eps)^{d-1}/d}\cdot \frac{K^i\eps}{d}\le 1.
$$
Thus, $n_{i-1}\le d^22^{(i-1)(d-1)}K^{-i}\eps^{-d}$.

Therefore, for $1\le i\le l$, the measure of $U_i$ is bounded above
by \bqn s_i\cdot (2^{-i}\eps)^d
&\le&  2^dn_{i-1}\cdot (2^{-i}\eps)^d\\
&\le& 2^d\cdot d^22^{(i-1)(d-1)}K^{-i}\eps^{-d}\cdot (2^{-i}\eps)^d\\
&=&2d^2(2K)^{-i}. \eqn For $i=0$, the measure of $U_0$ is trivially
bounded by 1.

Recall that for $0\le i<l$, $|\overline{f}-\underline{f}|\le
2K^{i+1}\eps$ on $U_i$. Also, on $U_l$, we have
$|\overline{f}-\underline{f}|\le 1$. Thus, \bq
\|\overline{f}-\underline{f}\|_p^p
&=&\int_{U_0}|\overline{f}-\underline{f}|^p
          +\sum_{i=1}^{l-1}\int_{U_i}|\overline{f}-\underline{f}|^p
          +\int_{U_l} |\overline{f}-\underline{f}|^p\nonumber\\
&\le& (2K \eps)^p+ \sum_{i=1}^{l-1}(2K^{i+1}\eps)^p \cdot 2d^2
          \left(2K\right)^{-i}+ 2d^2\left(2K\right)^{-l}\nonumber\\
&\le& (2K\eps)^p+
             2^{p+1}K^pd^2\sum_{i=1}^{l-1}\left(\frac{K^{p-1}}{2}\right)^i\eps^p
             +2d^2\left(2K\right)^{-l}.\label{f-f}
\eq

When $(d-1)p<d$, we have $d-1<\beta<\frac1{p-1}$. So,
$K=2^\beta<2^{1/(p-1)}$. Thus, $K^{p-1}/2<1$, and $\frac{1}{2K}\le
K^{-p}$. Therefore \bq \|\overline{f}-\underline{f}\|_p^p &\le& (2K
\eps)^p+ 2^{p+1}K^pd^2\cdot \frac{K^{p-1}}{2-K^{p-1}}\eps^p+2d^2
            \cdot K^{-pl}\nonumber\\
&\le& \left[(2K)^p+2^{p+1}K^pd^2\cdot
          \frac{K^{p-1}}{2-K^{p-1}}+2d^2\right]\eps^p\nonumber\\
&\le& c\eps^p\label{f-f1} \eq for some constant $c$ depending only
on $p$ and $d$, where in the second inequality we used the fact that
$K^{-l}\le \eps$.

When $(d-1)p>d$, we have $d-1>\beta>\frac1{p-1}$. So,
$K=2^\beta>2^{1/(p-1)}$, that is $K^{p-1}/2>1$. Hence, \bqn
\|\overline{f}-\underline{f}\|_p^p&\le& (2K \eps)^p+
2^{p+1}K^pd^2\cdot \frac{(K^{p-1}/2)^{l}}{K^{p-1}/2-1}\eps^p+2d^2
\cdot (2K)^{-l}\\
&\le& (2K \eps)^p+ \frac{2^{p+1}K^pd^2}{K^{p-1}/2-1}\cdot
K^{pl}\eps^p\cdot(2K)^{-l}+2d^2
\cdot (2K)^{-l}\\
&\le &(2K \eps)^p+c(2K)^{-l}\\  &\le&
(2K)^p\eps^p+c\eps^{1+1/\beta}\\
&\le& c'\eps^{1+1/\beta},\eqn for some constants $c,c'>0$ depending
only on $p$ and $d$, where in the third and fourth inequalities we
used the fact $1\le K^l\eps <K$ and in last inequality we used the
fact that $p>1+1/\beta$.

When $(d-1)p=d$, we have $K^{p-1}=2$, So, we obtain from (\ref{f-f})
that \bqn \|\overline{f}-\underline{f}\|_p^p &\le& (2K\eps)^p+
             2^{p+1}K^pd^2(l-1)\eps^p
             +2d^2\left(K^p\right)^{-l}\\&\le& c\eps^p\log 1/\eps,
\eqn for some constant $c>0$ depending only on $p$, where in the
last inequality we used the fact that $1\le K^{l}\eps<K$.\\

Summarizing, we obtain that \bq\|\overline{f}-\underline{f}\|_p\le
\threepiece{c\eps}{(d-1)p<d}{c\eps(\log
1/\eps)^{1/p}}{(d-1)p=d}{c\eps^{\frac{\beta+1}{p\beta}}}{(d-1)p>d}.\label{difference}\eq

\subsection{Bounds for $|\overline{\cS}|$ and $|\underline{\cS}|$}
We derive the upper bound for $|\overline{\cS}|$. The argument for bounding
$|\underline{\cS}|$ is almost identical.

Because all the selected cubes of side-length $\eps$ are chosen from
$n_0=\eps^{-d}$ cubes, there are no more than $2^{\eps^{-d}}$
different ways of selecting cubes of side-length $\eps$. For $1\le
i<l$, the selected cubes of side-length $2^{-i}\eps$ are chosen from
the $n_{i-1}$ cubes of side-length $2^{-i+1}\eps$ that were not
selected in the previous step, there are no more than
$2^{2^dn_{i-1}}$ different ways to select the cubes of side-length
$2^{-i}\eps$. Once the cubes are selected. For each $0\le i<l$, the
$s_i$ selected cubes of side-length $2^{-i}\eps$ can be grouped into
no more than $(2^{i}\eps^{-1})^{d-1}$ rows. Suppose row-$j$ contains
$r_j$ selected cubes. Because the values of $\overline{f}$ on these
$r_j$ cubes are in monotonic order, and are all chosen from $0$,
$K^i\eps$, $2K^i\eps$, ... $mK^i\eps$, where
$m=\floor{K^{-i}\eps^{-1}}$, the number of different ways of
assigning values of $\overline{f}$ on these $r_j$ cubes is bounded
by
$${r_j+\floor{K^{-i}\eps^{-1}}
\choose{\floor{K^{-i}\eps^{-1}}+1}} \le \max\{\exp(cr_j),
\exp(cK^{-i}\eps^{-1})\}<\exp(cr_j)\cdot \exp(cK^{-i}\eps^{-1}).$$
Thus, the number of different ways to assign the values of
$\overline{f}$ on the $s_i$ selected cubes of side-length
$2^{-i}\eps$ is bounded by \bqn
\prod_{j=1}^{(2^{i}\eps^{-1})^{d-1}}\left(\exp(cr_j)\cdot
\exp(cK^{-i}\eps^{-1})\right)&\le& \exp(cs_i)\cdot
\exp\left(c(2^{d-1}K^{-1})^i\eps^{-d}\right)\\
&\le& \exp\left(c'(2^{d-1}K^{-1})^i\eps^{-d}\right),\eqn where in
the inequality above, we used $s_i\le 2^dn_{i-1}$, and the estimate
$n_{i-1}\le d^22^{(i-1)(d-1)}K^{-i}\eps^{-d}$ obtained in $\S3.2$.

Hence, the total number of realizations of $\overline{f}$ is bounded
by \bq 2^{\eps^{-d}}e^{c'\eps^{-d}}
\prod_{i=1}^{l-1}\left[2^{2^dn_{i-1}}\cdot
\exp\left(c'(2^{d-1}K^{-1})^i\eps^{-d}\right)\right]\le
\exp\left(c''\sum_{i=0}^{l-1}(2^{d-1}K^{-1})^i\eps^{-d}\right),\label{S}\eq
where in the last inequality we again used the estimate $n_{i-1}\le
d^22^{(i-1)(d-1)}K^{-i}\eps^{-d}$.

When $(d-1)p>d$, $2^{d-1}>2^\beta=K$, we can bound the right hand
side of (\ref{S}) by \bqn
\exp\left(c'''[2^{d-1}/K]^l\eps^{-d}\right)\le
\exp\left(c'''\eps^{-(\beta+1)(d-1)/\beta}\right).\eqn

When $(d-1)p=d$, the upper bound of the right hand side of (\ref{S})
can be bounded by $ \exp\left(c'''\eps^{-d}\log 1/\eps\right)$.

When $(d-1)p<d$, $2^{d-1}/K<1$, and the upper bound of the right
hand of (\ref{S}) is bounded by $\exp(c'''\eps^{-d})$.

Summarizing, we obtain \bq \log |\overline{\cS}|\le
\threepiece{c'''\eps^{-d}}{(d-1)p<d} {c'''\eps^{-d}\log
1/\eps}{(d-1)p=d}{c'''\eps^{-(\beta+1)(d-1)/\beta}}{(d-1)p>d}.\label{card}
\eq

\subsection{Proof of Proposition 3.1}
Combining (\ref{difference}) and (\ref{card}), we have
$$
\log N_{[\,]}(\eps, \cF_d,\|\cdot\|_p )\le
\threepiece{c\eps^{-d}}{(d-1)p<d}{c\eps^{-d}(\log
1/\eps)^{1+d/p}}{(d-1)p=d}{c\eps^{-(d-1)p}}{(d-1)p>d},
$$
for all $\eps=2^{-n}$, $n\in \NN$. The monotonicity of bracketing
numbers implies that Proposition 3.1 holds for all $\eps<1$.
\\

\section{Critical Case}
We believe that the logarithmic factor in
Theorem~\ref{BracketingEntropyBoundTheorem} is not needed. In this
section, we prove that if we only consider the regular entropy, then
when $(d,p)\ne(2,2)$, the logarithmic factor can indeed be removed.
\begin{theorem} For $(d,p)\ne (2,2)$, there exist constants $c_1, c_2$ depending only on
$p$ and $d$ such that,
$$
{c_1\eps^{-\alpha}}\le \log N(\eps, \cF_d,\|\cdot\|_p)\le
{c_2\eps^{-\alpha}},
$$
where $\alpha=\max\{d,(d-1)p\}$.
\end{theorem}
\begin{proof}
In view of Theorem 1.1, it remains to show the upper bound for the
case $(d-1)p=d$, $d>2$. Let
$$
T=\{1_A: A=\{(x_1,x_2,...,x_d): f(x_1,x_2,...,x_d)\le \lambda\},
0\le \lambda\le 1, f\in \cF_d\}.
$$
Then clearly $\cF_d$ is the closed convex hull of $T$, that is
$\cF_d={\rm conv}(T)$.

For any $1_A\in T$, there exists $f\in \cF_d$, and $0\le \lambda\le
1$ such that
$$
A=\{(x_1,...,x_d): f(x_1,...,x_d)\le \lambda\}.
$$
By otherwise changing variable $t_i=1-x_i$, we can assume that $f$ is
non-decreasing with respect to every variable $x_i$, $1\le i\le d$.
Define $f_A$ on $[0,1]^{d-1}$ as follows:
$$
f_A(x_1,x_2, ... ,x_{d-1})=\twopiece{\max\{t: (x_1,...,x_{d-1},t)\in
A\}} {\mbox{if } \{t: (x_1,...,x_{d-1},t)\in A\}\ne
\emptyset}{0}{\mbox{if }\{t: (x_1,...,x_{d-1},t)\in A\}= \emptyset}.
$$
It is easy to check that $f_A\in \cF_{d-1}$. Furthermore, for all
$1_A, 1_B\in T$, $\|1_A-1_B\|_p=\|f_A-f_B\|_1^{1/p}$. Thus,
$$
N_{[\,]}(\eps, T,\|\cdot\|_p)=N_{[\,]}(\eps^p, \cF_{d-1},\|\cdot\|_1
).
$$
Therefore, by applying Proposition 3.1 for $\cF_{d-1}$ with $p=1$,
we have
$$
\log N(\eps, T,\|\cdot\|_p )\le \log N_{[\,]}(\eps, T,\|\cdot\|_p
)\le c\eps^{-(d-1)p}.
$$
Recall a general theorem of \cite{ckp:99} (see also \cite{carl:97})
that $$\log N(\eps, {\rm conv}(S) )=O(\eps^{-\sigma})$$ whenever
$\log N(\eps, S )=O(\eps^{-\sigma})$ for $\sigma>2$. Applying these
results we obtain
$$
\log N(\eps, \cF_d,\|\cdot\|_p )=\log N(\eps, {\rm
conv}(T),\|\cdot\|_p )\le c\eps^{-(d-1)p},
$$
for $(d-1)p=d>2$.
\end{proof}
When $(p,d)=(2,2)$, we have $(d-1)p=2$. It was proved in
\cite{gao:04} that
$$\log N(\eps, {\rm conv}(S))=O(\eps^{-2}(\log 1/\eps)^2)$$ whenever
$\log N(\eps,S)=O(\eps^{-2})$, and in general, this cannot be
improved. Note that this bound is exactly the bound we obtained
earlier using a direct construction. Thus, in this case, using convex
hulls does not improve the estimate.

\section{Rates of convergence for the Maximum Likelihood Estimator of a block decreasing density}

Biau and Devroye \cite{biau-devroye:03} showed that the minimax rate
of convergence for estimating a bounded block decreasing density
with $L_1$ risk is  $n^{1/(2+d)}$, and constructed histogram
estimators that attain this rate.   Here is a more precise
description of their result. Let ${\cal F}_B$ denote the class of
all block decreasing densities on the unit cube $[0,1]^d$ bounded by
$B$. Define the risk of the estimator $\widehat{f}_n$ when the true
density is $f \in {\cal B}$ by
\begin{eqnarray*}
R( \widehat{f}_n , f) = E_f  \left \{ \int_{\RR^d} | \widehat{f}_n (x) - f(x) | \, dx \right \} ,
\end{eqnarray*}
and the maximum (or ``worst case'') risk by
\begin{eqnarray*}
{\cal R} ( \widehat{f}_n , {\cal F}_B ) = \sup_{f \in {\cal F}_B} R( \widehat{f}_n , f) .
\end{eqnarray*}
The {\sl minimax risk} is  ${\cal R}_n ({\cal F}_B) = \inf_{\widehat{f}_n} {\cal R} ( \widehat{f}_n, {\cal F}_B)$.
\cite{biau-devroye:03} showed that for some constants $C_1$ and $C_2$,
\begin{eqnarray*}
{\cal R}_n ({\cal F}_B) \ge C_2 \left ( \frac{ C_1 S^d}{n} \right )^{1/(d+2)}
\end{eqnarray*}
where $S \equiv \log (1+B)$.  The resulting minimax lower bound rate
of convergence is $r_n^{mmlb} = n^{1/(2+d)} = n^{\gamma/(2
\gamma+1)} $ where $1/\gamma = d$. \cite{biau-devroye:03} also
constructed generalizations of the histogram estimators of Birge
\cite{birge:87b} which achieve this rate of convergence.

The MLE of a decreasing density on $[0,M]$ is well known to be
$n^{1/3}$ with respect to Hellinger and $L_1$ metrics:  see Birg\'e
\cite{birge:86}, \cite{birge:87a}, \cite{birge:89}. Although the MLE
of a block decreasing density has been initiated by Polonik
\cite{polonik:98}, the rate of convergence of the MLE in this
setting with respect to Hellinger or $L_1$ metrics is  apparently
unknown for $d\ge 2$. It is known from Birg\'e and Massart
\cite{birge-msrt:93} (see also \cite{vdv-w:96}, pages 326-327
together with
 Theorem 3.4.1, page 322) that
maximum likelihood estimators have a rate of convergence of at least
$r_n^{mle} = n^{\gamma/2}$ when the bracketing entropy with respect to the Hellinger metric $h$
 of the class of densities
${\cal P} $ satisfies
\begin{eqnarray}
\log N_{[\,]} (\epsilon, {\cal P} , h ) \le \frac{K}{\epsilon^{1/\gamma}}, \qquad   \ \epsilon > 0
\label{TypicalBracketingEntropyBound}
\end{eqnarray}
with $\gamma < 1/2$; here the Hellinger distance $h(P,Q)$ is given
by $h^2 (p,q) = \int [ \sqrt{p} - \sqrt{q}]^2 d\mu$ where $\mu $ is
any measure dominating both $P$ and $Q$ and $p$, $q$ are the
densities of $P, Q$ with respect to $\mu$.   From the results of
\cite{biau-devroye:03} it might be guessed that
(\ref{TypicalBracketingEntropyBound}) holds for ${\cal P} = {\cal
F}_B$ with $1/\gamma = d$, and this would lead to the rate of
convergence $r_n=n^{1/(2d)}$ for the MLE when $d \ge 2$. Our theorem
1.1 suggests that the rate of the convergence of the MLE (with
respect to Hellinger distance) is still slower than this for $d>2$, as is
shown in the following proposition.  We suppose that $X_1 , \ldots ,
X_n$ are i.i.d. $f \in {\cal F}_B$.

\begin{proposition}
\label{MLERateProp} Suppose that $\widehat{f}_n$ is the MLE of a
block decreasing density $f$ on $[0,1]^d$.  Then if $d\ge 3$
\begin{eqnarray}
n^{\frac{1}{4(d-1)}}  h ( \widehat{f}_n , f ) = O_p (1) .
\label{MLERateStatementDimensionsThreeAndUp}
\end{eqnarray}
If $d=2$, then
\begin{eqnarray}
\frac{n^{1/4} }{\log n }  h ( \widehat{f}_n , f ) = O_p (1) .
\label{MLERateStatementDimensionTwo} \label{MLERateStatement}
\end{eqnarray}
\end{proposition}


\begin{proof}
We use the results of Birg\'e and Massart \cite{birge-msrt:93} as
presented in section 3.4 of \cite{vdv-w:96}.  From Theorem 3.4.1,
page 322, with $\Theta_n$ taken to be
$$
{\cal P} = \{ p \ \mbox{a block-decreasing density on} \ [0,1]^d \ \mbox{bounded by} \ \ B \}
$$
it follows that we need to establish the inequalities of the first display of page 323.
These follow from Theorem 3.4.4, page 327, for the Hellinger distance $h$  by choosing
$p_n = p_0$ and taking ${\cal P}_n = {\cal P}$:  the resulting bound for
$E_{P_0} \| \GG_n \|_{{\cal M}_{\delta} } $ with
$$
{\cal M}_{ \delta} = \{ m_p = \log \frac{p+ p_0}{p_0} : \ p \in {\cal P} \}
$$
is of the form
\begin{eqnarray}
\tilde{J}_{[\, ]} ( \delta , {\cal P} , h) \left ( 1 + \frac{\tilde{J}_{[\, ]} ( \delta , {\cal P} , h)}{ \delta^2 \sqrt{n} } \right )
\equiv \phi_n (\delta)
\label{PhiOscillationBound}
\end{eqnarray}
where
$$
\tilde{J}_{[\, ]} ( \delta , {\cal P} , h) = \int_{c \delta^2}^{\delta} \sqrt{ 1 + \log N_{[\, ]} ( \epsilon, {\cal P}, h)} d \epsilon
$$
in view of the discussion on page 326 and \cite{birge-msrt:93},
Theorem 1, page 118. Since $\sqrt{p}$ is block-decreasing with bound
$\sqrt{B}$ if $p$ is block-decreasing with bound $B$, it follows
that
$$
\log N_{[\, ]} ( \epsilon , {\cal P}, h  ) = \log N_{[\, ]} ( \epsilon, {\cal P}^{1/2}, \| \cdot \|_2 )
=  \log N_{[\, ]} ( \epsilon/\sqrt{B}, {\cal P}^{1/2}/\sqrt{B} , \| \cdot \|_2 )
$$
where  $\| \cdot \|_2 $ is the $L_2 $ norm (with respect to Lebesgue
measure $\lambda$) and where ${\cal P}^{1/2}$ is the class of block
- decreasing functions with bound $\sqrt{B}$, and hence ${\cal
P}^{1/2}/\sqrt{B}$ is the class of block - decreasing functions with
bound $1$.  Thus for $d\ge 3$ we calculate, using Theorem 1.1 with
$p=2$,
\begin{eqnarray*}
\tilde{J}_{[\, ]} ( \delta , {\cal P} , h)
& = & \int_{c \delta^2}^{\delta} \sqrt{ 1 + \log N_{[\, ]} ( \epsilon, {\cal P}, h)} d \epsilon \\
& = & \int_{c \delta^2}^{\delta}
           \sqrt{ 1 + \log N_{[\, ]} ( \epsilon/\sqrt{B} , {\cal P}^{1/2}/\sqrt{B} , \| \cdot \|_2 )} d \epsilon \\
& \le  & \twopiece{\int_{c \delta^2}^{\delta}  \sqrt{ 1 +  c_2 B^{d-1} \epsilon^{ - 2(d-1)} } d \epsilon }{d>2}
{\int_{c \delta^2}^{\delta}  \sqrt{ 1 +  c_2 B \epsilon^{ - 2}(\log 1/\eps)^2 } d \epsilon }{d=2}\\
& \lesssim & \twopiece{\delta^{-2(d-2)}}{d>2}{(\log
1/\delta)^2}{d=2}
\end{eqnarray*}
where $f(x) \lesssim g(x)$ means $f(x) \le K g(x)$ for some constant $K$.
Plugging this into (\ref{PhiOscillationBound}) yields
$$
\phi_n (\delta ) = \delta^{-2(d-2)} \left ( 1 +
\frac{\delta^{-2(d-2)}}{\delta^2 \sqrt{n} } \right ) \,\,\,\, \mbox{
for }d>2,
$$
$$
\phi_n (\delta ) = (\log(1/\delta)^2 \left ( 1 +
\frac{(\log(1/\delta)^2}{\delta^2 \sqrt{n} } \right ) \,\,\,\,
\mbox{ for }d=2.
$$
It is easily verified that when $d>2$, $r_n^2 \phi_n ( 1/r_n )
\lesssim \sqrt{n}$ if $r_n = n^{\frac1{4(d-1)}}$. When $d=2$, $r_n^2
\phi_n ( 1/r_n ) \lesssim \sqrt{n}$ if $r_n = n^{\frac14}/\log n$.
Thus the rate of convergence of the MLE  is at least
$n^{\frac1{4(d-1)}}$ for $d>2$, and $n^{\frac14}/\log n$ for $d=2$.
\end{proof}

\bibliographystyle{ims}

\end{document}